\documentclass[12pt,twoside]{article} 
\usepackage{amsmath,amssymb,mathrsfs,graphicx}

\newtheorem{theorem}{Theorem}[section]
\newtheorem{lemma}[theorem]{Lemma}
\newtheorem{proposition}[theorem]{Proposition}

\newtheorem{definition}[theorem]{Definition}

\newtheorem{rmrk}[theorem]{Remark}

\DeclareMathAlphabet{\mathbfit}{OML}{cmm}{b}{it}

\makeatletter
\@addtoreset{equation}{section}
\makeatother

\newenvironment{remark}
{\begin{rmrk} \em}
{\end{rmrk}}

\newcommand{\fn} {function}
\newcommand{\me} {measure}
\newcommand{\tr} {trajector}
\newcommand{\erg} {ergodic}
\newcommand{\sy} {system}

\newcommand{\pr} {probability}
\newcommand{\dsy} {dynamical system}
\renewcommand{\o} {orbit}

\newcommand{\R} {\mathbb{R}}

\newcommand{\Z} {\mathbb{Z}}
\newcommand{\N} {\mathbb{N}}
\newcommand{\qed} {\hfill {\small Q.E.D.} \par\medskip}
\newcommand{\skippar} {\par\medskip}
\newcommand{\ds} {\displaystyle}
\newcommand{\proof} {\noindent \textsc{Proof.} }
\newcommand{\proofof}[1] {\noindent \textsc{Proof of {#1}.} }
\newcommand{\article}[3] {\textsc{{#1}}, {\itshape {#2}}, {{#3}}.}
\newcommand{\book}[3] {\textsc{{#1}}, {\itshape {#2}}, {{#3}}.}
\newcommand{\vol} {\textbf}
\newcommand{\eps} {\varepsilon}

\newcommand{\rset}[2] {\left\{ #1 \: \left| \: #2 \right. \! \right\} }

\newcommand{\into} {\longrightarrow}

\renewcommand{\emptyset} {\varnothing}

\newcommand{\ob} {observable}

\newcommand{\ps} {I}   % unit interval
\newcommand{\br} {\tau}   % branch fn
\newcommand{\pa} {\mathscr{P}}   % generic symbol for partitions
\newcommand{\ind} {\mathcal{J}}   % generic symbol for indices
\newcommand{\leb} {m}   % Lebesgue measure
\newcommand{\bj} {\mathbfit{j}}  % bold italic j
\newcommand{\bor} {\mathscr{B}}   % Borel sigma-algebra
\newcommand{\sca} {\mathscr{A}}   % generic sigma-algebra
\newcommand{\uc} {\mathbb{S}^1}   % unit circle
\newcommand{\bt} {\mathcal{T}}   % extended map

\begin{document}

\title{\textbf{A simple proof of the exactness \\
  of expanding maps of the interval \\
  with an indifferent fixed point}}

\author{
\scshape
Marco Lenci\,\thanks{
Dipartimento di Matematica, Universit\`a di Bologna,
Piazza di Porta San Donato 5, 40126 Bologna, Italy. 
E-mail: \texttt{marco.lenci@unibo.it}.}
\thanks{
Istituto Nazionale di Fisica Nucleare,
Sezione di Bologna, Via Irnerio 46,
40126 Bologna, Italy.}
}

\date{November 2015 \\[6pt] 
\normalsize{Final version to be published in \\
{\it Chaos Solitons \& Fractals}}}

\maketitle

\begin{abstract}
  Expanding maps with indifferent fixed points, a.k.a.\ intermittent
  maps, are popular models in nonlinear dynamics and infinite ergodic
  theory.  We present a simple proof of the exactness of a wide class
  of expanding maps of $[0,1]$, with countably many surjective
  branches and a strongly neutral fixed point in 0.
  
\end{abstract}

\section{Introduction}
\label{sec-intro}

Expanding maps with indifferent fixed points are very popular models
in nonlinear dynamics. Not only are they among the simplest chaotic
\dsy s whose physical \me\ may be infinite, they have also been used
to model anomalous transport in deterministic settings (see, e.g.,
\cite{gt}, \cite{gnz}, \cite[\S1.2.3.3]{bg}, \cite{zk}, \cite{k} and
references therein).

Constructing one such scheme is simple. Starting, say, with a
sufficiently regular map $T: \uc \into \uc$ with an indifferent fixed
point in $\bar{x}$, we lift $T$ to a map $\tau : [0,1) \into \R$ (that
is, after choosing an identification $\uc \cong [0,1)$, $\tau$ is
such that $T(x) = \tau(x)$ mod 1). Then we define $\bt : \R \into
\R$ as the translation invariant version of $\tau$ (namely, for $x
\in [k, k+1[$, $\bt(x) := \tau(x-k) +k$). Maps like $\bt$ are
called \emph{quasi-lifts} in \cite{l14}.

One is interested in the diffusive properties of the \tr ies of $\bt$,
for example as \fn s of an initial condition $x \in [0,1)$, chosen
w.r.t.\ the physical \me\ normalized to $[0,1)$.

Assume that $\tau$ was chosen so that $\tau(\bar{x}) = \bar{x}$.  Then
$\bt$ has countably many indifferent fixed points in $\bar{x}_k :=
\bar{x} + k$, $k \in \Z$. When $\bt^n(x)$ gets very close to one of
them, the \tr y will remain around it for a long while. After that, it
will reach a region where the modulus of the derivative of $\bt$ is
substantially different from 1, causing it to undergo an erratic, or
chaotic, motion. This will end when $\bt^{n'}(x)$, with $n' > n$, gets
again very close to a fixed point, and so on. So the typical \tr y
will alternate between almost constant stretches and random-looking
stretches.  In jargon, it will have an \emph{intermittent} behavior
(whence the name \emph{intermittent maps} for expanding maps with
indifferent fixed points). The statistical properties of
$(\bt^n(x))_n$, seen as a random process, may be very different from
the case of a similar map with expanding fixed points only: for
example, the scaling rate that is observed in the (generalized) CLT
and/or the mean square displacement might be $n^\gamma$, with $\gamma
\ne 1/2$. In this case one speaks of \emph{anomalous} diffusion
\cite{bg}.

This phenomenon can also be studied in terms of the original map $T$,
which we now view as a map $[0,1) \into [0,1)$ (via the aforementioned
identification $\uc \cong [0,1)$). Choose an initial condition $x \in 
[0,1)$ and let $\lfloor \cdot \rfloor$ denote the integer part of a
real number. The definitions of $\tau$ and $\bt$ easily imply that
$\bt(x) = \tau(x) = \lfloor \tau(x) \rfloor + T(x)$. Setting $f(x)
:= \lfloor \tau(x) \rfloor$ and iterating the procedure, we obtain
\begin{equation} \label{intro}
  \bt^n(x) = \sum_{j=0}^{n-1} f(T^j(x)) + T^n(x) ,
\end{equation}
Since $T^n(x)$ is bounded, the diffusive properties of $\bt^n(x)$ are
completely revealed by the Birkhoff sum of the observable $f: [0,1)
\into \Z$.

The latter may be called \emph{discrete displacement}, as it specifies
in what copy of the unit interval the dynamics is going to take place
at the next iteration.  One can also study Birkhoff sums of more
general \ob s $f$, taking values in $\Z$ or $\R$, regular or not
around $\bar{x}$, vanishing or not there. Each choice gives rises to
different statistical properties of the random varables $\sum_j f
\circ T^j$, which can always be viewed as the \tr ies of a given
extended \dsy\ (called a \emph{group extension} or \emph{skew product}
\cite[\S8.1]{a}).

\skippar

This preamble was meant to illustrate the importance, from the point
of view of applications, of the stochastic properties of interval maps
with indifferent fixed points. Exactness is one of the strongest of
these properties: for a non-singular \dsy, it means that the
\sy\ eventually loses all initial information---encoded in the form of
an absolutely continuous \pr\ \me\ for the choice of the initial
conditions. (The reader unfamiliar with this or other notions of the
theory of \dsy s is referred to the brief recapitulation of Section
\ref{app-basic} of the Appendix.)  If the \sy\ preserves a finite
absolutely continuous \me, exactness implies mixing of all orders
\cite{r,q}. If the relevant invariant \me\ is instead infinite,
exactness is perhaps the only notion of strong mixing whose definition
works well in infinite \erg\ theory too (see the discussion in
\cite{l13, l14}).

In this paper we deal with expanding Markov maps of the interval with
a finite number of indifferent fixed points. These maps are always
non-singular w.r.t.\ the Lebesgue \me\ and, in great generality,
possess a unique absolutely continuous invariant \me\ \cite{t80}.
Under some conditions on the nature of the fixed points, such \me\ is
infinite.  Famous examples are the Pomeau-Manneville maps \cite{pm}
and the Farey map (see, e.g., \cite{i} or \cite{ks}).

In 1983, Thaler proved that if a map as described above has surjective
branches only, then it is exact under very mild technical conditions
\cite{t83}.  This theorem is partly based on previous work by the same
author \cite{t80}. In his celebrated 1997 book, Aaronson extended the
result to a large class of Markov maps in a general setting
\cite[\S4]{a}.  Understandably, such general proof is rather
cumbersome. On the other hand, Thaler's original papers are not
straightforward either, as they involve non-standard types of induced
maps, a martingale convergence theorem and so on. (Recent proofs of
the exactness of specific maps, such as the Farey map \cite{ks} and
$\alpha$-Farey maps \cite{kms}, are not easily generalizable, or
especially simple either.)

The purpose of this note is to present a hands-on and relatively short
proof of the exactness of the simplest kinds of Markov maps of $[0,1)$
preserving an infinite \me: those defined by an indifferent fixed
point at 0 and a countable number of uniformly expanding surjective
branches.  This is not a serious restriction within Thaler's family,
as will become clear below.

What makes our key argument rather immediate is the use of a recent
criterion for exactness by Miernowski and Nogueira \cite{mn} (a
generalization of which we present in the Appendix). Understandably,
the argument needs distortion estimates. The ones we give here are
transparent---at least in this author's view---for they are based on a
simple estimate by Young \cite[\S6]{y}. To make this paper
self-contained, Young's proof is reported in the Appendix too.

\bigskip\noindent
\textbf{Acknowledgments.} I would like to thank Claudio Bonanno, Sara
Munday and Lai-Sang Young for useful discussions, and Roberto Artuso
for pointing out some relevant references. This work is part of my
activities within the \emph{Gruppo Nazionale di Fisica Matematica}
(INdAM, Italy). It was also partially supported by PRIN Grant
2012AZS52J\underline{ }001 (MIUR, Italy).

\section{Setup and result}
\label{sec-setup}

Many of the least common mathematical terms used in this section are
defined in Section \ref{app-basic} of the Appendix.

We assume that there is a partition $\pa := \{ I_j \}_{j \in \ind}$ of
$\ps := [0,1]$. The partition can be finite, in which case $\ind :=
\Z_N := \{ 0, 1, \ldots, N-1 \}$ or countable, in which case $\ind :=
\N$ (in our notation $0 \in \N$). The elements of the partition are
defined to be $I_j = [a_j, a_{j+1}]$, with $0 = a_0 < a_1 < \ldots <
a_k < \ldots$

Let $T : \ps \into \ps$ be a Markov map w.r.t.\ $\pa$, with the
following properties:
\begin{itemize} 
\item[(A1)] $T |_{(a_j, a_{j+1})}$ possesses an extension $\br_j : I_j
  \into I$ which is bijective and $C^2$ up to the boundary.
  
\item[(A2)] There exists $\Lambda > 1$ such that $| \br'_j (x) | \ge
  \Lambda$, for all $x \in I_j$ with $j \ge 1$.
  
\item[(A3)] There exists $K > 0$ such that $\ds \frac{ |\br''_j (x)| }
  { |\br'_j (x)|^2 } \le K$, for all $x$ and $j$.
  
\item[(A4)] $\br_0$ is convex with $\br_0(0) = 0$, $\br'_0(0) = 1$,
  $\br'_0(x) > 1$, for $x>0$, and $\br''_0(x) \sim x^{\beta}$, for $x
  \to 0^+$, for some $\beta \ge 0$.
\end{itemize}

It is proved in \cite{t80} that, under the above conditions, $T$
possesses an infinite invariant \me\ $\mu$, absolutely continuous
w.r.t.\ the Lebesgue \me\ $\leb$ and such that $h := d\mu / d\leb$ is
bounded on every $[\eps, 1]$. The arguments there, as well as in
\cite[\S1.5]{a}, prove that $\mu$ is unique up to factors.  In any
event, the point of view of this note is that the map $T$ is given
together with the \me\ $\mu$ it preserves, as is the case in many
applications. This way, none of our proofs depend on \cite{t80}.

\bigskip\noindent
\textbf{Terminology and conventions.} 
\begin{enumerate}
\item Unless it is important and clearly specified, neither our
  notation nor our language will mention null-\me\ sets. For example,
  we liberally say that $\pa$ is a partition of $\ps$ even though
  $I_j$ and $I_{j+1}$ intersect in a point; or we write $T I_j = \ps$
  even though this might be true only mod $\mu$.

\item Throughout the paper, the $\sigma$-algebra of reference for
  $\ps$ will be its Borel $\sigma$-algebra $\bor$. In fact, every time
  a $\sigma$-algebra is implied in the arguments, we shall always
  intend the Borel $\sigma$-algebra of the space at hand.
\end{enumerate}

This is our main result:

\begin{theorem} \label{thm-ex} 
  $T: I \into I$ is conservative and exact, w.r.t.\ $\mu$, or,
  equivalently, $\leb$.
\end{theorem}

\proof It is easy to check that $\forall x \in (0, a_1)$, $\exists n
\in \N$ such that $x < T(x) < \ldots < T^n(x) \not\in I_0$. So $J :=
\bigcup_{j \in \ind \setminus \{0 \}} I_j$ is a global cross-section,
in the sense that almost every \o\ of the \sy\ intersects it. Moreover
$\mu(J) < \infty$. Therefore, via the Poincar\'e Recurrence Theorem
applied to the map induced by $T$ on $J$, the \dsy\ is conservative.

As for the exactness, we are going to use the Miernowski-Nogueira
criterion \cite{mn}:

\begin{proposition} \label{prop-mn} 
  The non-singular and \erg\ \dsy\ $(X, \sca, \nu, T)$ is exact if,
  and only if, $\forall A \in \sca$ with $\nu(A) > 0$, $\exists n =
  n(A)$ such that $\nu( T^{n+1}A \cap T^n A ) > 0$.
\end{proposition}

A generalization of this criterion to the case of non-\erg\ maps is
given in Section \ref{app-mymn} of the Appendix.

\skippar

We need to define a more refined Markov partition for $T$. Let
$(b_k)_{k \in \N} \subset I_0$ be uniquely defined by $b_0 := a_1$ and
$T(b_{k+1}) = b_k$, with $b_{k+1} < b_k$. Now, for $k \in \Z^+$, set
$I_{-k} := [b_k, b_{k-1}]$. Then, $\pa_- := \{ I_j \}_{j \in \Z^-}$ is
a partition of $I_0$. So $\pa_o := \pa_- \cup \pa \setminus \{ I_0 \}$
is a partition of $\ps$. Its index set will be denoted $\ind_o := \Z^-
\cup \ind \setminus \{0\}$. $\pa_o$ is a Markov partition because $T(
I_{-1} ) = J$ and, for $k \ge 2$, $T( I_{-k} ) = I_{-k+1}$.

Let $\pa_o^n := \bigvee_{k=0}^{n-1} T^{-k} \pa_o$ denote the
refinement of $\pa_o$ relative to $T$ up to time $n$. For $\bj^n :=
(j_0, \ldots, j_{n-1}) \in (\ind_o)^n$, its elements are denoted
\begin{equation} \label{ex-5} 
  I_{\bj^n} := I_{j_0} \cap T^{-1} I_{j_1} \cap \cdots \cap T^{-n+1} 
  I_{j_{n-1}}
\end{equation}
(notice that there are many $\bj^n$ for which $I_{\bj^n} =
\emptyset)$.  Since $T$ is uniformly expanding away from 0, and since
a.a.\ \o s visit $J$ infinitely often, it is easily seen that, for any
sequence $(j_n)_{n \in \N} \subset \ind_o$,
\begin{equation} \label{ex-10}
  \lim_{n \to \infty} \leb ( I_{(j_0, \ldots, j_{n-1})} ) = 0.
\end{equation}

We now enter the core of the proof. Let $A$ be any positive-\me\ set.
Among the infinitely many density points of $A$, relative to $\leb$,
let us choose $x_0$ so that its \o\ intersects $J$ infinitely many
times (this is possible because $J$ is a global cross-section). Let
$(j_n)$ describe the itinerary of $x_0$ w.r.t.\ $\pa_o$, namely, $T^n
(x_0) \in I_{j_n}$, $\forall n \in \N$; equivalently, $x_0 \in
I_{(j_0, \ldots, j_{n-1})} = I_{\bj^n}$, $\forall n \in \N$. By
(\ref{ex-10}), using the notation of conditional \me,
\begin{equation} \label{ex-20}
  \lim_{n \to \infty} \leb ( A \,|\, I_{\bj^n} ) = 1.
\end{equation}
Moreover, we can assume that there exist $\bar{\jmath} \in \Z^+$ and a
subsequence $(j_{n_k})$ such $j_{n_k} = \bar{\jmath}$. In fact,
keeping in mind that $J = \bigcup_{j \in \ind \setminus \{0 \}} I_j$
is a global cross-section, if the \o\ of $x_0$ intersected each $I_j$,
with $j \ge 1$, only a finite number of times, then necessarily $T^n
(x_0) \to 1$, as $n \to \infty$. But $T$ is conservative and 1 is not
an atom of $\mu$, so there can only be a null-\me\ set of such points,
and we can pick a different $x_0$.

We need a distortion lemma, which will be proved in Section
\ref{sec-dist}.

\begin{lemma} \label{lem-dist}
  There exists $D > 1$ such that, for any $n \in \N$; any $\bj^{n+1} =
  (j_0, \ldots, j_n) \in (\ind_o)^{n+1}$ with $\leb( I_{\bj^{n+1}} ) >
  0$ and such that at least one of its components $j_k > 0$; and any
  $B \subseteq I_{\bj^{n+1}}$, one has:
  \begin{itemize}
  \item[(i)] $T^n B \subseteq I_{j_n}$;
  \item[(ii)] $\leb( T^n B \,|\, I_{j_n} ) \le D \, \leb( B \,|\,
    I_{\bj^{n+1}} )$;
  \item[(iii)] $\leb( T^{n+1} B ) \le D \, \leb( B \,|\, 
    I_{\bj^{n+1}} )$.
  \end{itemize}
\end{lemma}

From now till the end of the proof, to comply with one of the
hypotheses of the lemma, we always take $n \ge n_1$: that way, for $k
= n_1$, $j_k = \bar{\jmath} > 0$.

Applying Lemma \ref{lem-dist}\emph{(ii)} to $B := I_{\bj^{n+1}}
\setminus A$, whose conditional Lebsegue \me\ in $I_{\bj^{n+1}}$
vanishes by (\ref{ex-20}), and observing that $I_{j_n} \cap (T^n A)
\supseteq I_{j_n} \setminus T^n B$, we see that
\begin{equation} \label{ex-30}
  \lim_{n \to \infty} \leb ( T^n A \,|\, I_{j_n} ) = 1.
\end{equation}

Now we notice that $\exists \delta \in (0,1)$ such that, if $C \subset
I_{\bar{\jmath}}$ with $\leb( C \,|\, I_{\bar{\jmath}} ) > \delta$,
then $\leb( C \cap TC ) > 0$. (This is not hard to prove, using Lemma
\ref{lem-dist}\emph{(iii)} with $n=0$, $j_0 = \bar{\jmath}$, and $B =
I_{\bar{\jmath}} \setminus C$. The optimal estimate for $\delta$ is
found to be $D / ( D + \leb( I_{\bar{\jmath}} ) )$.)  Therefore,
choosing a sufficiently large $k$ such that, by (\ref{ex-30}), $\leb (
T^{n_k} A \,|\, I_{j_{n_k}} ) > \delta$, and since $j_{n_k} =
\bar{\jmath}$, we obtain $\leb(T^{n_k} A \cap T^{n_k + 1} A) > 0$,
which is the exactness condition of Proposition \ref{prop-mn}.

In order to apply that proposition, we still need to verify that $T$
is \erg. But this follows immediately from the above arguments. In
fact, if $A$ is an invariant set with $\leb(A)>0$, (\ref{ex-30}) gives
\begin{equation} \label{ex-40}
   \leb ( A \,|\, I_{\bar{\jmath}} ) = \lim_{k \to \infty} \leb 
   ( T^{n_k} A \,|\, I_{j_{n_k}} ) = 1.
\end{equation}
Using that $TA = A$, $T I_{\bar{\jmath}} = I$ and $T$ is non-singular,
(\ref{ex-40}) implies that $\leb(A) = \leb( TA \,|\, T
I_{\bar{\jmath}} ) = 1$.
\qed

\section{Distortion}
\label{sec-dist}

This section is entirely devoted to the proof of Lemma \ref{lem-dist}.
We will use standard techniques and variations thereof.

Firstly, \emph{(i)} follows from (\ref{ex-5}) since $B \subseteq
I_{\bj^{n+1}}$. Secondly, \emph{(ii)} comes from the following
distortion inequality: $\forall x,y \in I_{\bj^{n+1}}$,
\begin{equation} \label{dist-5}
  D^{-1} \le \frac{ |(T^n)' (x)| } { |(T^n)' (y)| } \le D.
\end{equation}
In fact,
\begin{equation} \label{dist-7}
\begin{split}
  \leb( T^n B \,|\, I_{j_n} ) &= \frac{\leb( T^n B )} 
  {\leb( T^n I_{\bj^{n+1}} )} = \frac{\int_B |(T^n)'| \, d\leb} 
  {\int_{I_{\bj^{n+1}}} |(T^n)'| \, d\leb} \\
  &\le \frac{ \max_B  |(T^n)'| \: \leb(B)} { \min_{I_{\bj^{n+1}}}  |(T^n)'| 
  \: \leb(I_{\bj^{n+1}}) } \le D \, \leb( B \,|\, I_{\bj^{n+1}} ) .
\end{split}
\end{equation}
Assertion \emph{(iii)} is derived in the same way from the inequality
\begin{equation} \label{dist-10}
  D^{-1} \le \frac{ |(T^{n+1})' (x)| } { |(T^{n+1})' (y)| } \le D,
\end{equation}
using that $T^{n+1} B \subseteq T^{n+1} I_{\bj^{n+1}} = I$.

Thus, we need to prove (\ref{dist-5}) and (\ref{dist-10}). We will
only write the proof of the latter, since the former is completely
analogous (and in fact implied by our proof, as will be clear).
Denoting for short $x_k := T^k(x)$ and $y_k := T^k(y)$, an easy
sufficient condition for (\ref{dist-10}) is
\begin{equation} \label{dist-20}
  \left| \sum_{k=0}^n \log \frac{ |T'(x_k)| } { |T'(y_k)| } \right| 
  \le C,
\end{equation}
where $C$ is a positive constant (whence $D := e^C$).

Observe that, by definition of $I_{\bj^{n+1}}$, the \o\ segments
$(x_k)_{k=0}^n$ and $(y_k)_{k=0}^n$ have the same itinerary
w.r.t\ $\pa_o$. We are going to parse them by grouping
\emph{excursions in $I_0$}, where we define an excursion in $I_0$ to
be an \o\ segment $\{ x_i, x_{i+1}, \ldots, x_j \}$ such that $x_i \in
J$ and $x_k \in I_0$, for all $i < k \le j$. The excursion is said to
be \emph{complete}, respectively \emph{partial}, if $x_{j+1} \in J$,
respectively $x_{j+1} \in I_0$.

Set $k_0 := 0$, and, recursively for $i \ge 1$, $k_i := \min \rset{k >
k_{i-1}} {x_k \in J}$ (the definition would be equivalent with $y_k$
in place of $x_k$).  This process stops when there are no more $k_i
\le n$ to define. We denote by $\ell$ the last index $i$ for which
$k_i$ has been defined, and also set $k_{\ell+1} := n+1$.

So each time frame $\{ k_i, k_i + 1, \ldots, k_{i+1} - 1\}$
corresponds to one of following four types of \o\ segments:

\begin{description}
\item[\sc Type 1:] The first segment of the parsing, which might not be 
an actual excursion in $I_0$, if $x_0, y_0 \not\in J$.

\item[\sc Type 2:] \emph{Bona fide} complete excursions in $I_0$, that is, 
complete excursions of cardinality bigger than 1. 

\item[\sc Type 3:] \emph{Degenerate} excursions, that is, single points 
in $J$ followed by points in $J$.

\item[\sc Type 4:] The last segment of the parsing, which might only 
be a partial excursion, if $x_{n+1}, y_{n+1} \not\in J$.
\end{description}

\begin{remark} \label{rk-nontrivial}
  The hypothesis that at least one of the $j_k$ is positive means
  that, for at least one $k$, $x_k, y_k \in J$. This implies that the
  parsing is not trivial, i.e., it cannot comprise just one segment.
  Otherwise, as will be clear below, certain estimates might be
  arbitrarily bad.
\end{remark}

We are going to show that there exist constants $\eta \in (0,1)$ and
$\kappa > 0$ such that, in each time frame of type 1-3, we have:
\begin{align}
  \label{est1}
  | x_{k_i} - y_{k_i} | &\le \eta | x_{k_{i+1}} - y_{k_{i+1}} | ; \\
  \label{est2} 
  \left| \sum_{k = k_i}^{k_{i+1} - 1} \log \frac{ |T'(x_k)| } { |T'(y_k)| } 
  \right| &\le \kappa | x_{k_{i+1}} - y_{k_{i+1}} | .
\end{align}
For the type 4 segment, we have
\begin{equation} \label{est3}
  \left| \sum_{k = k_\ell}^n \log \frac{ |T'(x_k)| } { |T'(y_k)| } \right| \le 
  \kappa.
\end{equation}

The estimates (\ref{est1})-(\ref{est3}) yield (\ref{dist-20}) because 
\begin{equation} \label{dist-50}
\begin{split}
  \left| \sum_{k=0}^n \log \frac{ |T'(x_k)| } { |T'(y_k)| } \right| &\le 
  \sum_{i=0}^{\ell-1} \left| \sum_{k = k_i}^{k_{i+1} - 1}  
  \log \frac{ |T'(x_k)| } { |T'(y_k)| } \right| + \left| \sum_{k = k_\ell}^n 
  \log \frac{ |T'(x_k)| } { |T'(y_k)| } \right| \\
  &\le \sum_{i=0}^{\ell-1} \kappa \eta^{\ell - i} \, | x_{k_{\ell+1}} -
  y_{k_{\ell+1}} | + \kappa \\
  &\le \frac{\kappa} {1 - \eta} =: C,
\end{split}
\end{equation}
where we have used that $| x_{k_{\ell+1}} - y_{k_{\ell+1}} | = |
x_{n+1} - y_{n+1} | \le 1$.

\skippar

Let us prove (\ref{est1})-(\ref{est2}) for each of the first three
types of \o\ segments, starting with the easiest.

\medskip\noindent 
\textsc{Type 3.} Since $x_{k_i}, y_{k_i} \in J$, (A2) yields
\begin{equation} \label{dist-60}
  | x_{k_i + 1} - y_{k_i +1} | \ge \Lambda  | x_{k_i} - y_{k_i} |.
\end{equation}  
Furthermore, let $j \in \ind$ be such that $x_{k_i},  y_{k_i}
\in I_{j}$. For some $\xi$ between $x_{k_i + 1}$ and $y_{k_i +1}$ 
one has
\begin{equation} \label{dist-70}
\begin{split}
  \left| \log \frac{ |T'(x_{k_i})| } { |T'(y_{k_i})| } \right| &=
  \left| \log \frac{ |T'( \br_{j}^{-1} (x_{k_i + 1}) )| } 
  { |T'( \br_{j}^{-1} (y_{k_i + 1}) )| } \right| \\ 
  &= \left| \frac{ T''( \br_{j}^{-1} (\xi) ) } { T'( \br_{j}^{-1} 
  (\xi) ) } \, \frac1 { T'( \br_{j}^{-1} (\xi) ) } \right| | x_{k_i + 1} - 
  y_{k_i +1} | \\
  &\le K \, | x_{k_i + 1} - y_{k_i +1} |,
\end{split}
\end{equation}  
by (A3). Since in this case $x_{k_i + 1} = x_{k_{i+1}}$ and 
$y_{k_i + 1}= y_{k_{i+1}}$, (\ref{est1})-(\ref{est2}) are shown.

\medskip\noindent 
\textsc{Type 2.} In this case too $x_{k_i}, y_{k_i} \in J$, therefore
(\ref{est1}) comes from (\ref{dist-60}) and the trivial inequality 
\begin{equation} \label{dist-75}
  | x_{k_{i+1}} - y_{k_{i+1}} | \ge | x_{k_i + 1} - y_{k_i +1} |. 
\end{equation} 

To show (\ref{est2}) we need the following lemma, which is practically
the same as \cite[\S6.2, Lem.~5]{y}. For the sake of completeness, we
give a proof in Section \ref{app-lsy} of the Appendix.

\begin{lemma} \label{lem-lsy}
  There exists $C' > 0$ such that, for all $j \ge 1$, $0 \le p \le j$,
  and $x,y \in I_{-j}$,
  \begin{displaymath}
    \left| \log \frac{ (T^p)' (x) } { (T^p)' (y) } \right| \le C'
    \frac{ |T^p(x) - T^p(y)| } { L_{p-j} } \le C',
  \end{displaymath}
  where, for $p \le j-1$, $L_{p-j} := |I_{p-j}| = b_{j-p-1} - b_{j-p}$
  and, for $p = j$, $L_0 := |J| = 1- a_1$ (observe that $T^p(x),
  T^p(y)$ belong to $I_{p-j}$ or $J$, respectively).
\end{lemma}

The l.h.s.\ of (\ref{est2}) can be estimated by
\begin{equation} \label{dist-80}
  \left| \sum_{k = k_i}^{k_{i+1} - 1} \log \frac{ |T'(x_k)| } 
  { |T'(y_k)| } \right| \le \left| \log \frac{ | T'(x_{k_i}) | } 
  { | T'(y_{k_i}) | } \right| + \left| \log \frac{ (T^p)' (x_{k_i +1}) } 
  { (T^p)' (y_{k_i +1}) } \right|,
\end{equation}  
with $p := k_{i+1} - k_i - 1$. (Notice that $(T^p)' (x_{k_i +1}) > 0$,
because $\{ x_{k_i +1}, \ldots, x_{k_{i+1} - 1} \} \subset I_0$, where
the map is increasing by (A4).) By (\ref{dist-70})-(\ref{dist-75}),
the first term of the r.h.s.\ of (\ref{dist-80}) is bounded above by
$K | x_{k_{i+1}} - y_{k_{i+1}} |$. For the second term we apply Lemma
\ref{lem-lsy}: in fact, $(T^p) (x_{k_i +1}) = x_{k_{i+1}}$, and the
same for $y_{k_i +1}$. Also, since $\{ x_{k_i}, \ldots, x_{k_{i+1} - 1} 
\}$ is a complete excursion, $x_{k_{i+1}}, y_{k_{i+1}} \in J$ by
construction: this means we apply the lemma in the case $p=j$.  The
second term in the assertion of Lemma \ref{lem-lsy} now reads $(C' /
|J|) | x_{k_{i+1}} - y_{k_{i+1}} |$. So (\ref{est2}) is proved for all
$\kappa \ge K + C' / |J|$.

\medskip\noindent 
\textsc{Type 1.} In this case, (\ref{est1}) is given by
\begin{equation} \label{dist-90}
   | x_{k_1} - y_{k_1} | \ge \bigg( \sup_{I_{-1}} |T'| \bigg) \,
   | x_{k_1 - 1} - y_{k_1 - 1} | \ge | x_0 - y_0 |.
\end{equation}  
As in the previous case, (\ref{est2}) follows from Lemma \ref{lem-lsy}, 
with $p = j := k_1$:
\begin{equation} \label{dist-100}
  \left| \sum_{k = 0}^{k_1 - 1} \log \frac{ |T'(x_k)| } { |T'(y_k)| } 
  \right| =  \left| \log \frac{ (T^{k_1})' (x_0) } { (T^{k_1})' (y_0) } 
  \right| \le  \frac{ C' } { |J| } | x_{k_1} - y_{k_1} |.
\end{equation}  

\medskip\noindent 
\textsc{Type 4.} It remains to verify (\ref{est3}) for the last
segment of the parsing. In analogy with (\ref{dist-80}),
\begin{equation} \label{dist-110}
  \left| \sum_{k = k_\ell}^n \log \frac{ |T'(x_k)| } { |T'(y_k)| } 
  \right| \le \left| \log \frac{ | T'(x_{k_\ell}) | } { | T'(y_{k_\ell}) 
  | } \right| + \left| \log \frac{ (T^p)' (x_{k_\ell +1}) } { (T^p)' 
  (y_{k_\ell +1}) } \right|,
\end{equation}  
with $p := k_{\ell+1} - k_\ell - 1 = n - k_\ell$. By (\ref{dist-70}),
the first term of the above r.h.s.\ is bounded above by $K | x_{k_\ell
+ 1} - y_{k_\ell + 1} | \le K$. The second term is bounded by $C'$
via the second inequality of Lemma \ref{lem-lsy}.

This concludes the proof of (\ref{est1})-(\ref{est3}) in all cases,
yielding (\ref{dist-50}), thus (\ref{dist-20}), thus Lemma
\ref{lem-dist}.

\appendix

\section{Appendix}

\subsection{Basic notions}
\label{app-basic}

We recall some basic notions of the mathematical theory of \dsy s that
have been used in the paper. Most of this material is presented, e.g.,
in \cite{a}.

\skippar

A {\bfseries \dsy} $(X, \sca, \nu, T)$ is defined by a \me\ space $(X,
\sca, \nu)$ and a map $T: X \into X$. We assume that $TX = X$.  $\sca$
is a $\sigma$-algebra defined on $X$ and $\nu$ is a $\sigma$-finite
\me\ for $(X,\sca)$ (this means that $\exists (A_n)_{n \in \N} \subset
\sca$, with $\nu(A_n) < \infty$ such that $\bigcup_n A_n = X$). The
\me\ of the space, $\nu(X)$, can be either finite or infinite; in the
former case, it is conventional to normalize $\nu$ so that $\nu(X) =
1$.  The map $T: X \into X$ is bi-measurable in the sense that,
$\forall A \in \sca$, both $T^{-1}A$ and $TA$ belong to $\sca$.

The \dsy, or the map, is called {\bfseries non-singular} if $A \in
\sca, \nu(A) = 0$ implies $\nu(T^{-1}A) = 0$. It is called {\bfseries
\me-preserving} if $\nu(T^{-1}A) = \nu(A)$, $\forall A \in \sca$
(equivalently, $\nu$ is said to be an {\bfseries invariant \me} for
$T$). Clearly, the latter property implies the former.

For a non-singular \dsy, a {\bfseries wandering set} is a measurable
set $W$ such that all the sets $\{ T^{-n}W \}_{n\in\N}$ are
disjoint. Points in $W$ have a non-recurrent behavior, insofar as, by
definition, $x \in W$ implies $T^n(x) \not\in W$, for all $n \ge
0$. It is always possible to partition $X$ into two parts
$\mathcal{D}$ and $\mathcal{C}$, defined up to null-\me\ sets, such
that every wandering set is contained in $\mathcal{D}$ (mod $\nu$).
$\mathcal{D}$ is called the {\bfseries dissipative part} of $X$ and
can be always be represented as a countable disjoint union of
wandering sets. $\mathcal{C}$ is called the {\bfseries conservative
part} of $X$ and it is where the recurrent behavior takes place. By
definition, in fact, every $A \subseteq \mathcal{C}$ is recurrent in
the sense of Poincar\'e, i.e., almost every $x \in A$ is such that
$T^n(x) \in A$, at a countable number of times $n$. A \dsy\ is called
{\bfseries conservative} if $X = \mathcal{C}$ and {\bfseries
dissipative} if $X = \mathcal{D}$. A finite-\me-preserving \sy\ is
always conservative (Poincar\'e Recurrence Theorem).

The \dsy\ is called {\bfseries \erg} if all the invariant sets are
trivial, namely, $T^{-1}A = A$ mod $\nu$ implies that either $\nu(A) =
0$ or $\nu(X \setminus A) = 0$.  This is equivalent to saying that the
{\bfseries invariant $\sigma$-algebra}
\begin{equation} \label{basic-inv}
    \mathscr{I} := \rset{A \in \sca} {T^{-1}A = A \ \mathrm{mod}\ \nu}
\end{equation}
is trivial (i.e., it contains only zero-\me\ sets and their
complements). Observe that, in the infinite-\me\ case, this is a
stronger notion than the classical definition whereby the time (i.e.,
Birkhoff) average of any \ob\ is constant almost everywhere.  For
example, in the case where $\nu$ is an infinite invariant \me, the
fact that
\begin{equation} \label{basic-birk-inf}
    \lim_{n \to \infty} \frac1n \sum_{k=0}^{n-1} f(T^k(x)) = 0,
\end{equation}
for any integrable \fn\ $f: X \into \R$ and almost every $x \in X$
(depending on $f$), does not imply that $T$ is \erg.

The {\bfseries tail $\sigma$-algebra} of a non-singular \dsy\
is defined to be
\begin{equation} \label{basic-tail}
    \mathscr{T} := \bigcap_{n=0}^\infty T^{-n} \sca.
\end{equation}
If $T: X \into X$ is a bijection, clearly $\mathscr{T} = \sca$, so
this quantity is only relevant for non-invertible maps. In rough
terms, we might say that the structure of $\mathscr{T}$ represents the
order that persists when the dynamics evolves and chaos is
produced. Equivalently, the information that we obtain by observing
the \sy\ via $\mathscr{T}$-measurable \fn s is the information that
comes from the infinite past of the dynamics and is not increased by
further observations in time. Regardless, it is easy to see that
$\mathscr{I} \subseteq \mathscr{T}$. A non-singular \dsy\ is called
{\bfseries exact} if $\mathscr{T}$ is trivial. Hence, exactness
implies \erg ity.

If $T$ is non-singular and $X$ admits a finite or countable partition
$\pa$ such that, for any element $E \in \pa$, $TE$ is a union of
elements of $\pa$ and the restriction $T|_E$ is invertible mod $\nu$,
then $T$ is called a {\bfseries Markov map} relative to the {\bfseries
  Markov partition} $\pa$. (Depending on the context and the author,
more technical conditions are required in the definition of Markov
map. The overly general definition that we give here is sufficient for
our illustrative purposes.)

The most common examples of Markov maps are those defined on an
interval, with a Markov partition made up of sub-intervals in whose
interior $T$ is smooth. Such are the \sy s discussed in this paper.  A
map of this kind is called {\bfseries expanding} if $|T'(x)| > 1$, for
all $x$ where $T'(x)$ is defined.  It is called {\bfseries uniformly
expanding} if $\exists \lambda>1$ such that $|T'(x)| \ge \lambda$
for all $x$ as above. A {\bfseries fixed point} $\bar{x}$ (i.e.,
$T(\bar{x}) = \bar{x}$) is called {\bfseries expanding}, respectively
{\bfseries contracting}, if $|T'(\bar{x})| > 1$, respectively $<
1$. It is called {\bfseries indifferent}, or {\bfseries neutral}, if
$|T'(\bar{x})| = 1$. It is sometimes said that a fixed point is
strongly neutral if $T''(x)$ is regular around $\bar{x}$.

\subsection{A criterion for the exact components}
\label{app-mymn}

In this section we generalize Proposition \ref{prop-mn} to the case of
non-\erg\ maps, obtaining a characterization of the exact components
of non-singular maps.

\begin{definition} \label{def-ai}
  Let $(X, \sca, \nu, T)$ be a non-singular \dsy\ on a $\sigma$-finite
  \me\ space (cf.\ Section \ref{app-basic}).  We say that $A \in
  \sca$, with $\nu(A) > 0$, is {\bfseries asymptotically intersecting}
  w.r.t.\ the given \dsy\ if $\exists n = n(A)$ such that $\nu(
  T^{n+1}A \cap T^n A ) > 0$.  By the non-singularity of $T$, this is
  equivalent to $\nu( T^{k+1}A \cap T^k A ) > 0$, for all $k \ge n$.
\end{definition}

\begin{proposition} \label{prop-mymn} 
  In the framework of Definition \ref{def-ai}, let $\mathscr{I}$ and 
  $\mathscr{T}$ denote, respectively, the invariant and the tail 
  $\sigma$-algebras (cf.\ (\ref{basic-inv}), (\ref{basic-tail})). The 
  following holds:
  \begin{itemize}
  \item[(i)] if every positive-\me\ $A \in \mathscr{T}$ is
    asymptotically intersecting, then $\mathscr{I} = \mathscr{T}$;
    
  \item[(ii)] if $\mathscr{I} = \mathscr{T}$, then every
    positive-\me\ $A \in \sca$ is asymptotically intersecting.
  \end{itemize}
\end{proposition}

\begin{remark}
  Observe that \emph{(i)} is a stronger statement than the converse of
  \emph{(ii)}: in particular, combining \emph{(i)} and \emph{(ii)}, we
  see that if every set in the tail $\sigma$-algebra is asymptotically
  intersecting, so is every measurable set. Also, using the fact that
  any power of an exact map is exact, and vice-versa, it is easy to
  show that if $\mathscr{I} = \mathscr{T}$ then, $\forall A \in \sca$,
  with $\nu(A) > 0$, and $\forall \ell \in \Z^+$, $\exists n = n(A,
  \ell)$ such that $\nu( T^{k+j}A \cap T^k A ) > 0$, whenever $k \ge
  n$ and $1 \le j \le \ell$. In any event, Proposition \ref{prop-mn}
  is now a corollary of Proposition \ref{prop-mymn}.
\end{remark}

\proofof{Proposition \ref{prop-mymn}} We remark that the techniques
used here come entirely from \cite[Lem.~2.1]{mn}.

The proof of \emph{(i)} is already contained in \cite{l14}. We report
it here for the sake of completeness.  Take $B \in \mathscr{T}$. We
set out to prove that $B \in \mathscr{I}$.  If $\nu(B) = 0$, then
$\nu(T^{-1}B) = 0$ and $B \in \mathscr{I}$. So we assume that $\nu(B)
> 0$. It is a known simple fact that, for all $B \in \mathscr{T}$ and
$k \in \N$,
\begin{equation} \label{mn-10}
  B = T^{-k} \, T^k B.
\end{equation}
We want to show that $B = TB$ mod $\nu$. This and (\ref{mn-10}) will
imply that $T^{-1} B = T^{-1} TB = B$ mod $\nu$, whence $B \in
\mathscr{I}$, as desired.

Set $A := B \setminus TB \in \mathscr{T}$.  By (\ref{mn-10}), for all
$n \ge 0$,
\begin{equation} \label{mn-20}
  A = T^{-n} \, T^n B \setminus T( T^{-n-1} \, T^{n+1} B) =  T^{-n} 
  (T^n B \setminus T^{n+1} B), 
\end{equation}
whence
\begin{equation} \label{mn-30}
  T^n A = T^n B \setminus T^{n+1} B.
\end{equation}
Applying (\ref{mn-30}) with $n+1$ in lieu of $n$ implies that $T^{n+1}
A \subseteq T^{n+1} B$, which, compared again to (\ref{mn-30}), gives
$T^{n+1} A \cap T^n A = \emptyset$. Since this holds for all $n \in
\N$, the hypotheses imply that $\nu(A) = 0$. Thus, $B \subseteq TB$
mod $\nu$.

Analogously, setting $A' := TB \setminus B$, we get that, for all $n
\ge 0$, $T^n A' = T^{n+1} B \setminus T^n B$, whence $T^n A' \subseteq
T^{n+1} B$ and $T^{n+1} A' = T^{n+2} B \setminus T^{n+1}
B$. Therefore, $T^{n+1} A' \cap T^n A' = \emptyset$.  For the same
reasons as before, $TB \subseteq B$ mod $\nu$, which completes the
proof of assertion \emph{(i)}.

\skippar

As for \emph{(ii)}, assume by contradiction that there exists $A \in
\sca$ with $\nu(A) > 0$ such that $\nu( T^{n+1}A \cap T^n A ) = 0$,
for all $n \in \N$. We want to show that this is incompatible with
$\mathscr{I} = \mathscr{T}$.

Since $T$ is non-singular, the above assumption implies that
\begin{equation} \label{mn-105}
  \nu( T^{-n}  T^{n+1}A \cap A ) \le \nu( T^{-n} T^{n+1}A \cap 
  T^{-n} T^n A ) = 0.
\end{equation}
Therefore, setting
\begin{equation} \label{mn-120}
  B := \bigcup_{n \in \N} T^{-n} T^{n+1}A,
\end{equation}
we have
\begin{equation} \label{mn-110}
  \nu( B \cap A ) = 0.
\end{equation}

The sequence of sets in the r.h.s.\ of (\ref{mn-120}) is increasing,
so $B = \bigcup_{n \ge k} T^{-n} T^{n+1}A \in T^{-k} \sca$, for all $k
\in \N$, whence $B \in \mathscr{T}$.  On the other hand, if
$\mathscr{T} = \mathscr{I}$,
\begin{equation} \label{mn-130}
  B = T^{-1} B = \bigcup_{n \in \N} T^{-n-1} T^{n+1}A \supseteq A,
\end{equation}
which contradicts (\ref{mn-110}) because $\nu(A) > 0$.
\qed

\subsection{Young's distortion estimate}
\label{app-lsy}

Here we prove Lemma \ref{lem-lsy}, copying almost verbatim the proof
of Lemma 5, \S6.2 in \cite{y}.

\bigskip\noindent 
\textbf{Terminology.}\ In what follows, we write $f_n \sim g_n$ to
mean that $\exists \kappa_2 > \kappa_1 > 0$ such that $\kappa_1 < |
f_n/g_n | < \kappa_2$, for all $n$ (possibly with some restrictions,
if so specified); the same goes for other integer indices, such as
$k$, or $i$. Also, we write $f(x) \sim g(x)$ to mean that $\kappa_1 <
| f(x)/g(x) | < \kappa_2$ holds for all $x \in (0, a_1)$ (in all the
cases below, this will be equivalent to the asymptotics $x \to 0^+$).
\bigskip

Set $\alpha := 1/(\beta + 1)$ and, for $n \ge 1$, $\Delta n^{-\alpha}
:= n^{-\alpha} - (n+1)^{-\alpha}$. Observe that
\begin{equation} \label{lsy-5}
  \Delta n^{-\alpha} \sim n^{-(\alpha + 1)} =( n^{-\alpha} 
  )^{(\alpha+1) / \alpha}. 
\end{equation}
For $k \in \N$, let $n_k$ be the unique index such that
\begin{equation} \label{lsy-10}
  b_k \in [(n_k + 1)^{-\alpha}, n_k^{-\alpha}).
\end{equation}
By (A4), $T(x) - x \sim x^{\beta+2}$. This, the definition of $b_k$, 
and (\ref{lsy-5})-(\ref{lsy-10}) imply that
\begin{equation} \label{lsy-20}
\begin{split}
  \Delta b_k &:= b_{k-1} - b_k = T(b_k) - b_k \\
  &\sim b_k^{\beta+2} \sim ( n_k^{-\alpha} )^{\beta+2} \\
  &\sim ( \Delta n_k^{-\alpha} )^{(\beta+2) \alpha / (\alpha + 1)} \\
  &= \Delta n_k^{-\alpha}.
\end{split}
\end{equation}

Since, for every fixed positive integer $l$, $\Delta (n+l)^{-\alpha}
\sim \Delta n^{-\alpha}$ and $\Delta b_{k+l} \sim \Delta b_k$
(respectively, as \fn s of $n$ and $k$), the above shows that each
$I_{-k} = [b_k, b_{k-1}]$ intersects at most a bounded number of
intervals $[(n + 1)^{-\alpha}, n^{-\alpha})$, and vice-versa.

Recall that $x,y \in I_{-j}$. For $0 \le i \le p-1$, there exists
$\xi_i$ between $T^i(x)$ and $T^i(y)$ (hence $\xi_i \in I_{i-j}$) such
that
\begin{equation} \label{lsy-25}
  \log T'( T^i(x) ) - \log T'( T^i(y) ) = \frac{ T''(\xi_i) } 
  {T'(\xi_i) } \left( T^i(x) - T^i(y) \right).
\end{equation}
But, for all $i \in \{ 0, \ldots, p-1 \}$, $T''(\xi_i) \sim
\xi_i^\beta \sim b_{j-i}^\beta$; $T'(\xi_i) \sim 1$; and $| T^i(x) -
T^i(y)| \le L_{i-j} \sim \Delta b_{j-i} \sim b_{j-i}^{\beta+2}$. All
this implies that, for any $0 \le q \le p$,
\begin{equation} \label{lsy-30}
\begin{split}
  \left| \log \frac{ (T^q)' (x) } { (T^q)' (y) } \right| &\le 
  \sum_{i=0}^{q-1} \frac{ T''(\xi_i) } {T'(\xi_i) } \left| T^i(x) 
  - T^i(y) \right| \\
  &\le C_1 \sum_{i=0}^{q-1} b_{j-i}^{2\beta+2} 
  \le C_2 \sum_{i=0}^{q-1} n_{j-i}^{-\alpha (2\beta+2)} \\
  &\le C_2 \sum_{n=1}^\infty n^{-2} := C_3,
\end{split}
\end{equation}
where we have used the considerations of the first part of the proof
and $C_1, C_2$ are suitable positive constants.

The distortion inequality (\ref{lsy-30}) holds for a generic pair $x,y
\in I_{-j}$, not necessarily the one given in the statement of Lemma
\ref{lem-lsy}. By standard arguments---as in
(\ref{dist-5})-(\ref{dist-20})---it gives
\begin{equation} \label{lsy-38}
  e^{-C_3} \, \frac{ |x -y| } {L_{-j}} \le \frac{ | T^q(x) - T^q(y) | } 
  {L_{q-j}} \le e^{C_3} \, \frac{ |x -y| } {L_{-j}}.
\end{equation}
Comparing the above expression for a generic $q = i \in \{ 0, \ldots,
p-1 \}$ with the same for $q = p$, we see that, for all $0 \le i \le
p-1$,
\begin{equation} \label{lsy-40}
  \frac{ | T^i(x) - T^i(y) | } {L_{i-j}} \sim \frac{ | T^p(x) - 
  T^p(y) | } {L_{p-j}}.
\end{equation}

Using (\ref{lsy-40}) in the first line of (\ref{lsy-30}), with $q =
p$, together with some of the above estimates, yields
\begin{equation} \label{lsy-50}
\begin{split}
  \left| \log \frac{ (T^p)' (x) } { (T^p)' (y) } \right| &\le C' 
  \sum_{i=0}^{p-1} b_{j-i}^\beta \, L_{i-j} \, \frac{ | T^p(x) - 
  T^p(y) | } { L_{p-j} } \\
  &\le C' \frac{ | T^p(x) - T^p(y) | } { L_{p-j} }, 
\end{split}
\end{equation}
where $C'$ is another positive constant, and the last inequality
follows from $b_k^\beta < 1$ and $\sum_{k \in \N} L_{-k} = 1$.
\qed

\end{document}